\documentclass[11pt]{article}

\usepackage{amsmath}
\usepackage{amsfonts}
\usepackage{amssymb}

\newtheorem{theorem}{Theorem}
\newtheorem{corollary}{Corollary}
\newtheorem{definition}{Definition}

\newtheorem{lemma}{Lemma}
\newtheorem{proposition}{Proposition}
\newtheorem{remark}{Remark}
\newenvironment{proof}[1][Proof]{\noindent\textbf{#1.} }{\ \rule{0.5em}{0.5em}}

\def\l{\mathcal{L}_{q,1,v}}
\def\L{\mathcal{L}_{q,2,v}}

\def\F{\mathcal{F}_{q,v}}
\def\I{\infty}

\begin{document}

\title{\bf Heisenberg Uncertainty Principle for the $q$-Bessel Fourier transform}
\date{ }
\author{Lazhar Dhaouadi \thanks{%
Institut Pr\'eparatoire aux Etudes d'Ing\'enieur de Bizerte
(Universit\'e du 7 novembre à Carthage). Route Menzel Abderrahmene
Bizerte, 7021 Zarzouna, Tunisia. \quad\quad\quad\quad\quad\quad
E-mail lazhardhaouadi@yahoo.fr}} \maketitle

\begin{abstract}
In this paper we uses an I.I. Hirschman-W. Beckner entropy argument
to give an uncertainty inequality for the $q$-Bessel Fourier
transform:
$$
\mathcal{F}_{q,v }f(x)=c_{q,v }\int_{0}^{\infty }f(t)j_{v
}(xt,q^{2})t^{2v +1}d_{q}t,
$$
where $j_v(x,q)$ is the normalized Hahn-Exton $q$-Bessel function.
\end{abstract}

\section{Introduction}

I.I.Hirschman-W. Beckner entropy argument is one further variant of
Heisenberg's uncertainty principle.

Let $\widehat{f}$ be the Fourier transform of $f$ defined by
$$
\widehat{f}(x)=\int_{-\I}^{\I} f(y)e^{2i\pi xy}f(y)dy,\quad
x\in\mathbb R.
$$
If $f\in L^2(\mathbb R)$ with $L^2$-norme $\|f\|_2=1$, then by
Plancherel's theorem $\|\widehat{f}\|_2=1$, so that $|f(x)|^2$ and
$|\widehat{f}(x)|^2$ are probability frequency functions. The
variance of a probability frequency $g$ is defined by
$$
V[g]=\int_{\mathbb R} x^2g(x)dx.
$$
The Heisenberg uncertainty principle can be stated as follows
\begin{equation}
V[|f|^2]V[|\widehat{f}|^2]\geq\frac{1}{16\pi^2}.
\end{equation}
If $g$ is a probability frequency function, then the entropy of $g$
is defined by
$$
E(g)=\int_{\mathbb R} g(x)\log(x)dx.
$$
With $f$ as above, Hirschman [10] proved that
\begin{equation}
E(|f|^2)+E(|\widehat{f}|^2)\leq 0.
\end{equation}
By an inequality of Shannon and Weaver it follows that (2) implies
(1). Using the Babenko-Beckner inequality
$$
\|\widehat{f}\|_{p'}\leq A(p)\|f\|_p,\quad  1<p<2,\quad
A(p)=\left[p^{1/p}(p')^{-1/p'}\right]^{1/2},
$$
in Hirschman's proof of (2) another uncertainty inequality is
deduced. For more detail the reader can consult [8,10,11].

In this paper we use I.I. Hirschman entropy argument de give an
uncertainty inequality for the $q$-Bessel Fourier transform (also
called $q$-Hankel transform).

Note that other versions of the Heisenberg uncertainty principle for
the $q$-Fourier transform have recently appeared in the literature
[1,2,6]. There are some differences of the results cited above and
our result:

$\bullet$ In [1] the uncertainty inequality is established  for the
$q$-cosine and $q$-sine transform but here is established for the
$q$-Bessel transform.

$\bullet$ In [2] the uncertainty inequality is for the $q^2$-Fourier
transform but here is for the $q$-Hankel transform.

$\bullet$ In [6] the uncertainty inequality is established for
functions in $q$-Schwartz space. In this paper the uncertainty
inequality is established for functions in $\L$ space.

The inequality discuss here is a quantitative uncertainty principles
which give an information about how a function and its $q$-Bessel
Fourier transform relate. A qualitative uncertainty principles give
an information about how a function (and its Fourier transform)
behave under certain circumstances. A classical qualitative
uncertainty principle called Hardy's theorem. In [4,7] a $q$-version
of the Hardy's theorem for the $q$-Bessel Fourier transform was
established.

In the end, our objective is to develop a coherent harmonic analysis
attached to the $q$-Bessel operator
$$
\Delta _{q,v}f(x)=\frac{1}{x^{2}}\left[
f(q^{-1}x)-(1+q^{2v})f(x)+q^{2v}f(qx)\right].
$$
Thus, this paper is an opportunity to implement the arguments  of
the $q$-Bessel Fourier analysis proved before, as the Plancherel
formula, the positivity of the $q$-translation operator, the
$q$-convolution product, the $q$-Gauss kernel...

\section{The $q$-Bessel Fourier  transform}

\bigskip In the following we will always assume $0<q<1$ and $v>-1$. We denote by%
\begin{equation*}
\mathbb{R}_{q}=\left\{ \pm q^{n},\text{ \ }n\in \mathbb{Z}\right\}
,\text{ \ \ \ }\mathbb{R}_{q}^{+}=\left\{ q^{n},\text{ \ }n\in
\mathbb{Z}\right\} .
\end{equation*}%
For more informations on the $q$-series theory the reader can see
the references [9,12,14] and the references [3,5,13] about the
$q$-bessel Fourier analysis. Also for details of the proofs of the
following results in this section can be fond in [3].

\begin{definition}
\bigskip The q-Bessel operator is defined as follows
\begin{equation*}
\Delta _{q,v}f(x)=\frac{1}{x^{2}}\left[ f(q^{-1}x)-(1+q^{2v})f(x)+q^{2v}f(qx)%
\right].
\end{equation*}
\end{definition}

\begin{definition}The normalized q-Bessel function of
Hahn-Exton is defined by
$$
j_{v}(x,q^{2})=\sum_{n=0}^{\infty }(-1)^{n}\frac{q^{n(n+1)}}{
(q^{2v+2},q^{2})_{n}(q^{2},q^{2})_{n}}x^{2n}.
$$
\end{definition}

\begin{proposition}The function
$$
x\mapsto j_{v}(\lambda x,q^{2})
$$
is the eigenfunction of the operator $\Delta_{q,v}$ associated with
the eigenvalue $ -\lambda ^{2}$.
\end{proposition}

\begin{definition}
The q-Jackson integral of a function $f$ defined on $\mathbb{R}_{q}$
is
\begin{equation*}
\int_{0}^{\infty }f(t)d_{q}t=(1-q)\sum_{n\in
\mathbb{Z}}q^{n}f(q^{n}).
\end{equation*}
\end{definition}

\begin{definition}
We denote by $\mathcal L_{q,p,v}$  the space of even functions $f$
defined on $\mathbb{R}_q$ such that
$$
 \|f\|_{q,p,v}=\left[\int_0^{\infty}|f(x)|^px^{2v+1}d_qx\right]^{1/p}<\infty.
$$
\end{definition}

\begin{definition}
We denote by ${\mathcal C}_{q,0}$ the space of even functions
defined on $\mathbb R_q$ tending to $0$ as $x\rightarrow\pm\infty$
and continuous at $0$ equipped with the topology of uniform
convergence. The space ${\mathcal C}_{q,0}$ is complete with respect
to the norm
$$
\|f\|_{q,\infty}=\sup_{x\in\mathbb R_q}|f(x)|.
$$
\end{definition}

\begin{definition}
The q-Bessel Fourier transform $\mathcal{F}_{q,v}$ (also called
q-Hankel transform) is defined by
\begin{equation*}
\mathcal{F}_{q,v}f(x)=c_{q,v}\int_{0}^{\infty
}f(t)j_{v}(xt,q^{2})t^{2v+1}d_{q}t,\text{ \ \ \ }\forall x\in
\mathbb{R}_{q}.
\end{equation*}
where
$$
c_{q,v} =
\frac{1}{1-q}\frac{(q^{2v+2};q^2)_\infty}{(q^{2};q^2)_\infty}.
$$
\end{definition}

\begin{proposition}
Let $f\in \mathcal{L}_{q,1,v}$ then $\mathcal{F}_{q,v}f$ existe and $%
\mathcal{F}_{q,v}f\in \mathcal{C}_{q,0}$.
\end{proposition}

\begin{definition}
The $q$-translation operator is given as follows
\begin{equation*}
T_{q,x}^{v}f(y)=c_{q,v}\int_{0}^{\infty }\mathcal{F}%
_{q,v}f(t)j_{v}(yt,q^{2})j_{v}(xt,q^{2})t^{2v+1}d_{q}t\quad \forall
f\in\l .
\end{equation*}
\end{definition}

\begin{definition}
The operator $T_{q,x}^{v}$ is said  positive if $T_{q,x}^{v}f\geq 0$
when $ f\geq 0$ for all $x\in \mathbb{R}_{q}.$ We denote by $Q_{v}$
the domain of positivity of $T_{q,x}^{v}$
\begin{equation*}
Q_{v}=\left\{ q\in ]0,1[,\text{ \ \ \ \ \ }T_{q,x}^{v}\text{ is
positive} \right\} .
\end{equation*}
\end{definition}
In the following we assume that $q\in Q_{v}.$

\begin{proposition}
If $f\in \mathcal{L}_{q,1,v}$ then%
\begin{equation*}
\int_{0}^{\infty }T_{q,x}^{v}f(y)y^{2v+1}d_{q}y=\int_{0}^{\infty
}f(y)y^{2v+1}d_{q}y.
\end{equation*}
\end{proposition}

\begin{definition}
The $q$-convolution product is defined as follows
\begin{equation*}
f\ast _{q}g(x)=c_{q,v}\int_{0}^{\infty
}T_{q,x}^{v}f(y)g(y)y^{2v+1}d_{q}y.
\end{equation*}
\end{definition}

\begin{proposition}
Let $f,g\in \mathcal{L}_{q,1,v}$ then $f\ast _{q}g\in\l$  and we have%
\begin{equation*}
\F(f\ast _{q}g)=\F(g)\times \F(f).
\end{equation*}
\end{proposition}

\begin{proposition}
Let $f\in \mathcal{L}_{q,1,v}$ and $g\in \mathcal{L}_{q,2,v}$ then
$f\ast _{q}g\in\L$  and we have
\begin{equation*}
\mathcal{F}_{q,v}(f\ast _{q}g)=\mathcal{F}_{q,v}(f)\times
\mathcal{F}_{q,v}(g).
\end{equation*}
\end{proposition}

\begin{theorem}
The $q$-Bessel Fourier transform $\F$ satisfies

\item 1. $\F$ sends $\L$ to $\L$.

\item 2. For $f\in\mathcal L_{q,2,v}$, we have
$$
\|\F(f)\|_{q,2,v}=\|f\|_{q,2,v}.
$$

\item 3. The operator $\F:\L\rightarrow\L$ is bijective and
$$
\F^{-1}=\F.
$$
\end{theorem}

\bigskip

\begin{proposition} Given $1<p\leq 2$ and
$\frac{1}{p}+\frac{1}{\overline{p}}=1$. If $f\in\mathcal L_{q,p,v}$
then
$$
\mathcal F_{q,v}(f)\in\mathcal L_{\overline{p},2,v}
$$
and
$$
\|\mathcal F_{q,v}(f)\|_{q,\overline{p},v}\leq
B_{q,v}^{(\frac{2}{p}-1)}\|f\|_{q,p,v},
$$
where
$$
B_{q,v}=\frac{1}{1-q}\frac{(-q^2;q^2)_\infty(-q^{2v+2};q^2)_\infty}{(q^2;q^2)_\infty}.
$$
\end{proposition}

\begin{definition}
The $q$-exponential  function is defined by
$$
 e(z,q)=\sum_{n=0}^\infty\frac{z^n}{(q,q)_n}=\frac{1}{(z;q)_\infty},\quad|z|<1.
$$
\end{definition}

\begin{proposition}The $q$-Gauss kernel
$$
G^v(x,t^2,q^2)=\frac{(-q^{2v+2}t^2,-q^{-2v}/t^2;q^2)_\infty}{(-t^2,-q^2/t^2;q^2)_\infty}e\Big(-\frac{q^{-2v}}{t^2}x^2,q^2\Big),
\quad\forall x,t\in\mathbb R_q^+
$$
satisfies
$$
\mathcal F_{q,v}\left\{e(-t^2y^2,q^2)\right\}(x)=G^v(x,t^2,q^2),
$$
and for all function $f\in\mathcal L_{q,2,v}$
$$
\lim_{n\rightarrow \infty}\|G^v(x,q^{2n},q^2)*_q f-f\|_{q,2,v}=0.
$$
\end{proposition}

\section{Uncertainty Principle}

The following Lemma are  crucial for the proof of our main result.
First we enunciate the Jensen's inequality

\begin{lemma}
Let $\gamma$ be a probability measure on $\mathbb R_q^+$. Let $g$ be
a convex function on a subset  $I$ of $\mathbb{R}$. If $\psi:\mathbb
R_q^+\rightarrow I$ satisfies
$$
\int_0^\infty \psi(u)d\gamma(u)\in I,
$$
then we have
$$
g\left(\int_0^\infty \psi(x) d\gamma(x)\right)\leq\int_0^\infty
g\circ \psi(x)d\gamma(x).
$$
\end{lemma}

\begin{proof}Let
$$
t=\int_0^\infty \psi(u)d\gamma(u).
$$
There exist $c\in\mathbb R$ such that for all $y\in I$ it holds
$$
g(y)\geq g(t)+c(y-t).
$$
Now let $y=\psi(x)$ we obtain
$$
g\left( \psi(x)\right)\geq g(t)+c(\psi(x)-t).
$$
Integrating both sides and using the special value of t gives
$$
\int_0^\infty g\left ( \psi(x)\right)d\gamma(x)\geq \int_0^\infty
[g(t)+c(\psi(x)-t)]d\gamma(x)=g(t).
$$
This finish the proof.
\end{proof}

\begin{lemma}Let $f$ be an even function defined on $\mathbb R_q$.
Assume $\psi:\mathbb R\rightarrow \mathbb R_+$ is a convexe function
and $\psi\circ f\in\l$. If $\varrho_n$ is a sequence of non-negative
function in $\l$ such that
$$
\F(\varrho_n)(0)=c_{q,v}\int_0^\I \varrho_n(x)x^{2v+1}d_qx=1
$$
and $\varrho_n*_qf\rightarrow f$ then
$\psi\circ\Big(\varrho_n*_qf\Big)$ is in $\l$ and
$$
\lim_{n\rightarrow\I}\int_0^\I
\psi\circ\Big(\varrho_n*_qf\Big)(x)x^{2v+1}d_qx=\int_0^\I\psi\circ
f(x)x^{2v+1}d_qx.
$$
\end{lemma}

\begin{proof}For a given $x$ and by Proposition 3 we have
$$
c_{q,v}\int_0^\I T^v_{q,x}\varrho_n(y)y^{2v+1}d_qy=1
$$
From the positivity of $T^v_{q,x}$ we see that
$$
c_{q,v}T^v_{q,x}\varrho_n(y)y^{2v+1}d_qy
$$
is a probability measure on $\mathbb R_q^+$. The following holds by
Jensen's Inequality
$$\aligned
\psi\circ\Big(\varrho_n*_qf\Big)(x)&=\psi\left[c_{q,v}\int_0^\I
f(y)T^v_{q,x}\varrho_n(y)y^{2v+1}d_qy\right]\\
&\leq c_{q,v}\int_0^\I\psi\circ
f(y)T^v_{q,x}\varrho_n(y)y^{2v+1}d_qy\\
&=\varrho_n*_q\psi\circ f(x).
\endaligned$$
By the use of the Fatou's Lemma and Proposition 4 we obtain
$$\aligned
&\int_0^\I\psi\circ f(x)x^{2v+1}d_qx\\
&=\int_0^\I
\liminf_{n\rightarrow\I}\psi\circ\Big(\varrho_n*_qf\Big)(x)x^{2v+1}d_qx\\
&\leq
\liminf_{n\rightarrow\I}\int_0^\I\psi\circ\Big(\varrho_n*_qf\Big)(x)x^{2v+1}d_qx\\
&\leq
\limsup_{n\rightarrow\I}\int_0^\I\psi\circ\Big(\varrho_n*_qf\Big)(x)x^{2v+1}d_qx\\
&\leq
\lim_{n\rightarrow\I}\int_0^\I\varrho_n*_q\psi\circ f(x)x^{2v+1}d_qx\\
&=
\frac{1}{c_{q,v}}\lim_{n\rightarrow\I}\F(\varrho_n)(0)\times\F\Big(\psi\circ f\Big)(0)\\
&=\int_0^\I\psi\circ f(x)x^{2v+1}d_qx.
\endaligned$$
This finish the proof.
\end{proof}

\begin{definition}For a positive function $\phi$  define the entropy of $\phi$
to be
$$
E(\phi)=\int_{0}^{\infty}\phi(x)\log\phi(x)x^{2v+1}d_{q}x.
$$
$E(\phi)$ can have values in $[-\I,\I]$.
\end{definition}

\begin{remark}For a given $c\in\mathbb{R}_{q}^{+}$ let
$$
d\gamma(x)=k_{c}^{-1}\exp\left( -\left\vert cx\right\vert
^{a}\right) x^{2v+1}d_{q}x
$$
where
$$
\sigma_a=\int_{0}^{\infty}\exp\left( -\left\vert x\right\vert
^{a}\right) x^{2v+1}d_{q}x,\quad k_{c}=\frac{\sigma_a}{c^{2v+2}}.
$$
Then $d\gamma(x)$ is a probability measure on $\mathbb R_q^+$.
\end{remark}

\begin{lemma}Let $a>0$. For a positive function $\phi\in\mathcal
L_{q,1,v}$ such that $$\|\phi\|_{q,1,v}=1$$ and
$$
M_a(\phi)=\left( \int_{0}^{\infty}\left\vert x\right\vert
^{a}\phi(x)x^{2v+1}d_{q}x\right) ^{\frac{1}{a}}
$$
is finite,  we have
\begin{equation}
-E(\phi )\leq \log k_{c}+c^{a}M_{a}^{a}(\phi ).
\end{equation}
\end{lemma}

\begin{proof}
Indeed, defining
$$
\psi(x)=k_{c}\exp\left( \left\vert cx\right\vert ^{a}\right)
\phi(x),
$$
From Remark 1 we see that
$$
\int_{0}^{\infty}\psi(x)d\gamma(x)=1.
$$
According to the fact that $g:t\mapsto t\log t$ is convex on
$\mathbb{R}_+^{*},$ so Jensen's inequality gives
$$
g\left[\int_{0}^{\infty}\psi(x)d\gamma(x)\right]\leq\int_{0}^{\infty}g\circ\psi(x)d\gamma(x).
$$
Hence,
$$
0=\left[\int_{0}^{\infty}\psi(x)d\gamma(x)\right]\log\left[ \int_{0}^{\infty}\psi(x)d%
\gamma(x)\right] \leq\int_{0}^{\infty}\psi(x)\log\psi(x)d\gamma(x).
$$
This implies
\begin{eqnarray*}
0 &\leq &\int_{0}^{\infty }\phi (x)\log \left[ k_{c}\exp \left( \left\vert
cx\right\vert ^{a}\right) \phi (x)\right] x^{2v+1}d_{q}x \\
&=&\int_{0}^{\infty }\phi (x)\left[ \log k_{c}+\left\vert
cx\right\vert ^{a}+\log \phi (x)\right] x^{2v+1}d_{q}x.
\end{eqnarray*}

$$
0\leq \log k_{c}+c^{a}\int_{0}^{\infty }\left\vert x\right\vert
^{a}\phi (x)x^{2v+1}d_{q}x+\int_{0}^{\infty }\phi (x)\log \phi
(x)x^{2v+1}d_{q}x.
$$
In the end
$$
0\leq \log k_{c}+c^{a}M_{a}^{a}(\phi )+E(\phi ).
$$
This finish the proof.
\end{proof}

\begin{lemma}Let $f\in\l\cap\L$ then we have
\begin{equation}
E\left( \left\vert f\right\vert ^{2}\right) +E\left( \left\vert
\mathcal{F} _{q,v}f\right\vert ^{2}\right) \leq2\|f\|^2_{q,v,2}\log
\Big(B_{q,v}\|f\|^2_{q,v,2}\Big).
\end{equation}
\end{lemma}

\begin{proof}
H\"{o}lder inequality implies that $f$ will be in
$\mathcal{L}_{q,p,v}$ for $1<p\leq 2.$ With
$$
\frac{1}{p}+\frac{1}{\overline{p}}=1,
$$
Hausdorff-Young's inequality (Proposition 6) tells us that
$\mathcal{F}_{q,v}f$ \ is in $\mathcal{L}_{q,\overline{p},v}$. So we
can define the functions
$$
A(p)=\int_{0}^{\infty }\left\vert f(x)\right\vert ^{p}d_qx\text{ \ \ and }B(%
\overline{p})=\int_{0}^{\infty }\left\vert \mathcal{F}_{q,v}f(x)\right\vert
^{\overline{p}}x^{2v+1}d_{q}x\text{.}
$$
Now define
\begin{eqnarray*}
C(p) &=&\log \left\Vert\mathcal{F}_{q,v}f\right\Vert
_{q,\overline{p},v}-\log
\left( B_{q,v}^{\frac{2}{p}-1}\left\Vert f\right\Vert _{q,p,v}\right) \\
&=&\frac{1}{\overline{p}}\log B(\overline{p})-\frac{1}{p}\log
A(p)-\left( \frac{2}{p}-1\right) \log B_{q,v}.
\end{eqnarray*}
By Hausdorff-Young's inequality
$$
C(p)\leq 0,\text{ for }1<p<2,
$$
and by Plancherel equality (Theorem 1 part 2)
$$
C(2)=0.
$$
Then
$$
C'(2^-)\geq0.
$$
On the other hand for $1<p<2$ we have
$$
C^{\prime}(p)=\frac{\overline{p}^{\prime}}{\overline{p}}\frac{B^{\prime
}(
\overline{p})}{B(\overline{p})}-\frac{\overline{p}^{\prime}}{\overline
{p} ^{2}}\log B(\overline{p})-\frac{1}{p}\frac{A^{\prime}(p)}{A(p)}+
\frac{1}{p^{2}}\log A(p)+\frac{2}{p^{2}}\log B_{q,v}.
$$
The derivative of $\overline{p}$ with respect to $p$ is
$$
\overline{p}^{\prime}=-\frac{1}{(p-1)^{2}}.
$$
For a given $x>0$ we have
$$
\lim_{p\rightarrow2}\frac{x^{p}-x^{2}}{p-2}=x^{2}\log x.
$$
Then
$$\aligned
A^{\prime}(2^{-})&=\lim_{p\rightarrow2^{-}}\frac{A(p)-A(2)}{p-2}
=\frac{1}{2}E\left( \left\vert f\right\vert ^{2}\right),
\endaligned$$

$$
B^{\prime}(2^{+})=\lim_{\overline{p}\rightarrow2^{+}}\frac{B(\overline {p}%
)-B(2)}{\overline{p}-2}=\frac{1}{2}E\left( \left\vert
\mathcal{F}_{q,v}f\right\vert ^{2}\right) .
$$
Since
$$p\mapsto\frac{x^{p}-x^{2}}{p-2}
$$
is an increasing function, the exchange of the signs limit and
integral is valid sense. On the other hand
$$
\lim_{p\rightarrow2^{-}}A(p)=\|f\|^2_{q,v,2},
\quad\lim_{\overline{p}\rightarrow2^{+}}B(p)=\|\mathcal{F}_{q,v}f\|^2_{q,v,2}=\|f\|^2_{q,v,2}.
$$
So
$$
C^{\prime}(2^{-})=\lim_{p\rightarrow2^{-}}\frac{C(p)-C(2)}{p-2}=-\frac{1}{2\|f\|^2_{q,v,2}}%
\left[ A^{\prime}(2^{-})+B^{\prime}(2^{+})\right]
+\frac{1}{2}\log\Big(B_{q,v}\|f\|^2_{q,v,2}\Big).
$$
Therefore
$$
A^{\prime}(2^{-})+B^{\prime}(2^{+})-\|f\|^2_{q,v,2}\log
\Big(B_{q,v}\|f\|^2_{q,v,2}\Big)\leq0,
$$
and then
$$
E\left( \left\vert f\right\vert ^{2}\right) +E\left( \left\vert
\mathcal{F} _{q,v}f\right\vert ^{2}\right) \leq2\|f\|^2_{q,v,2}\log
\Big(B_{q,v}\|f\|^2_{q,v,2}\Big).
$$
This finish the proof.
\end{proof}

\begin{lemma}
Let $f\in\L$ then we have
\begin{equation}
E\left( \left\vert f\right\vert ^{2}\right) +E\left( \left\vert
\mathcal{F} _{q,v}f\right\vert ^{2}\right) \leq2\|f\|^2_{q,v,2}\log
\Big(B_{q,v}\|f\|^2_{q,v,2}\Big).
\end{equation}
\end{lemma}

\begin{proof}Assume that $E(|f|^2)$ and $E(|\F f|^2)$ are defined
and then approximate $f$ by functions in $\mathcal{L}_{q,1,v}\cap
\mathcal{L}_{q,2,v}.$ Let
$$
h_{n}(x)=e(-q^{2n}x^{2},q^{2}).
$$
The function $h_n$ is in $\L$ then $h_nf\in\l$. On the other hand
$h_n\in\mathcal{C}_{q,0}$ then $h_nf\in\mathcal{L}_{q,2,v}$. We
obtain
$$
h_nf\in\mathcal{L}_{q,1,v}\cap\mathcal{L}_{q,2,v}.
$$
The following holds by $(2)$
\begin{equation}
E\left(\left\vert h_nf\right\vert ^{2}\right) +E\left( \left\vert
\mathcal{F}_{q,v}(h_nf)\right\vert ^{2}\right)
\leq2\|h_nf\|^2_{q,2,v}\log \Big(B_{q,v}\|h_nf\|^2_{q,2,v}\Big).
\end{equation}
One can see by the Lebesgue Dominated Convergence Theorem that
\begin{equation}
\lim_{n\rightarrow \infty}\|h_nf\|_{q,2,v}=\|f\|_{q,2,v}
\end{equation}
and
\begin{equation}
\lim_{n\rightarrow \infty }E\left( \left\vert h_{n}f\right\vert
^{2}\right) =E\left( \left\vert f\right\vert ^{2}\right).
\end{equation}
By the use of Proposition 5 and the inversion formula (Theorem 1
part 3) we see that
$$
\F(h_nf)=\F h_{n}\ast _{q}\F f.
$$
We will prove that
$$
\lim_{n\rightarrow \infty }E\left( \left\vert \mathcal{F}_{q,v}h_{n}\ast _{q}%
\mathcal{F}_{q,v}f\right\vert ^{2}\right) =E\left( \left\vert \mathcal{F}%
_{q,v}f\right\vert ^{2}\right) .
$$
The functions
$$
\phi _{1}(x) =x^2\log ^{+}\left\vert x\right\vert \text{ and }\phi
_{2}\left( x\right) =x^2\left( -\log ^{-}\left\vert x\right\vert
+\frac{3}{2}\right),
$$
are convex on $\mathbb{R}$, where
$$
\log ^{+}x\text{ }=\text{max}\left\{ 0,\log x\right\} \text{ and }\log ^{-}x%
\text{ }=\text{min}\left\{ 0,\log x\right\} .
$$
Note that
$$
2\phi_1(x)-2\phi_2(x)+3x^2=x^2\log|x|^2.
$$
Since

$\bullet$ From the inversion formula we see that
$$c_{q,v}\int_{0}^{\infty
}\mathcal{F}_{q,v}h_{n}(t)t^{2v+1}d_{q}t=h_n(0)=1.$$

$\bullet$ The function $\F h_{n}\geq 0$.

$\bullet$ The functions $\phi_i$ are convex on $\mathbb{R}$.

$\bullet$ $E(\F f)$ is finite then $\phi_i(\F f)$ is in $\l$.

$\bullet$ From Proposition 7 we have
$$
\lim_{n\rightarrow\I}\F
h_{n}*_q\F f(x)=\F f(x)
$$
we deduce that $\F h_n$ and $\phi_i$ satisfy the conditions of Lemma
2. Then we obtain
$$
\lim_{n\rightarrow\I}\int_0^\I \phi_i\circ(\F h_n*_q\F
f)(x)x^{2v+1}d_qx=\int_0^\I \phi_i\circ(\F f)(x)x^{2v+1}d_qx,\quad
i=1,2.
$$
It also hold
$$\aligned
E\left( \left\vert \mathcal{F}_{q,v}f\right\vert ^{2}\right)
&=2\int_{0}^{\infty }\phi_{1}\left( \mathcal{F}_{q,v}f\right)
x^{2v+1}d_{q}\\
&-2\int_{0}^{\infty }\phi _{2}\left( \mathcal{F}_{q,v}f\right)
x^{2v+1}d_{q}x+3\left\Vert \mathcal{F}_{q,v}f\right\Vert
_{q,2,v}^{2},
\endaligned$$
and
$$\aligned
E\left( \left\vert \mathcal{F}_{q,v}h_{n}\ast _{q}\mathcal{F}%
_{q,v}f\right\vert ^{2}\right)&=2\int_{0}^{\infty }\phi _{1}\left(
\mathcal{F}_{q,v}h_{n}\ast _{q}\mathcal{F}_{q,v}f\right) x^{2v+1}d_{q}x \\
&-2\int_{0}^{\infty }\phi _{2}\left( \mathcal{F}_{q,v}h_{n}\ast _{q}%
\mathcal{F}_{q,v}f\right) x^{2v+1}d_{q}x \\
&+3\left\Vert \mathcal{F}_{q,v}h_{n}\ast
_{q}\mathcal{F}_{q,v}f\right\Vert _{q,2,v}^{2}.
\endaligned$$
Then
\begin{equation}
\lim_{n\rightarrow \infty }E\left( \left\vert \mathcal{F}_{q,v}h_{n}\ast _{q}%
\mathcal{F}_{q,v}f\right\vert ^{2}\right) =E\left( \left\vert \mathcal{F}%
_{q,v}f\right\vert ^{2}\right) .
\end{equation}
With $(6)$ and the limits $(7)$, $(8)$ and $(9)$ we complete the
proof of $(5)$.

\bigskip
Note that these limits also hold in the case where $E(|f|^2)$ and
$E(|\F f|^2)$ are $\I$ or $-\I$.
\end{proof}

\bigskip
Now we are in  position to state and prove the uncertainty
inequality for the $q$-Bessel Fourier transform.

\begin{theorem}Given $a, b > 0$. Then for all $c,d\in\mathbb R_q^+$ satisfying
$$
0<B^2_{q,v}\frac{\sigma_a\sigma_b}{(cd)^{2v+2}}<1,
$$
the following hold for any function $f\in\mathcal L_{q,2,v}$
$$
c^{a}\left\Vert x^{a/2}f\right\Vert _{q,2,v}^{2}+d^{b}\left\Vert
x^{b/2} \mathcal{F}_{q,v}f\right\Vert _{q,2,v}^{2}\geq -\log \left(
B^2_{q,v}\frac{\sigma_a\sigma_b}{(cd)^{2v+2}}\right) \left\Vert
f\right\Vert _{q,2,v}^{2}.
$$
\end{theorem}

\begin{proof}
Assume that $\|f\|_{q,2,v}=1$. By $(3)$ we can write
$$
-E(\left\vert f\right\vert ^{2})\leq\log k_{c}+c^{a}\left\Vert
x^{a/2}f\right\Vert _{q,2,v}^{2}
$$

$$
-E\left( \left\vert \mathcal{F}_{q,v}f\right\vert ^{2}\right) \leq
\log k_{d}+d^{b}\left\Vert x^{b/2}\mathcal{F}_{q,v}f\right\Vert
_{q,2,v}^{2}.
$$
Which implies with $(5)$
\begin{eqnarray*}
-2\log B_{q,v} &\leq &-E\left( \left\vert f\right\vert ^{2}\right)
-E\left(
\left\vert \mathcal{F}_{q,v}f\right\vert ^{2}\right) \\
&\leq &\log \Big(k_{c}k_{d}\Big)+c^{a}\left\Vert x^{a/2}f\right\Vert
_{q,2,v}^{2}+d^{b}\left\Vert x^{b/2}\mathcal{F}_{q,v}f\right\Vert
_{q,2,v}^{2}.
\end{eqnarray*}
By replacing $f$ \ by $\dfrac{f}{\left\Vert f\right\Vert _{q,2,v}}$ we get%
$$
c^{a}\left\Vert x^{a/2}f\right\Vert _{q,2,v}^{2}+d^{b}\left\Vert x^{b/2}%
\mathcal{F}_{q,v}f\right\Vert _{q,2,v}^{2}\geq -\log \left(
B^2_{q,v}k_{c}k_{d}\right) \left\Vert f\right\Vert _{q,2,v}^{2}.
$$
This finish the proof.
\end{proof}

\begin{corollary}There exist $k>0$ such that for any function
$f\in\mathcal L_{q,2,v}$ we have
$$
\|xf\|_{q,2,v}\|x\mathcal F_{q,v}f\|_{q,2,v}\geq k\|f\|^2_{q,2,v}.
$$
\end{corollary}

\begin{proof}Let $a=b=2$ and $c=d$ then by Theorem 3
$$
\left\Vert xf\right\Vert _{q,2,v}^{2}+\left\Vert x
\mathcal{F}_{q,v}f\right\Vert _{q,2,v}^{2}\geq -\frac{1}{c^2}\log
\left( B^2_{q,v}\frac{\sigma_2^2}{c^{4(v+1)}}\right) \left\Vert
f\right\Vert _{q,2,v}^{2},
$$
where
$$
0<\left( B^2_{q,v}\frac{\sigma_2^2}{c^{4(v+1)}}\right)<1.
$$
Now put
$$
f_t(x)=f(tx),\quad t\in\mathbb R_q^+,
$$
then
$$
\mathcal F_{q,v}f_t(x)=\frac{1}{t^{2v+2}}\mathcal
F_{q,v}f(x/t),\quad \|x\mathcal
F_{q,v}f_t\|^2_{q,2,v}=\frac{1}{t^{2v}}\|\mathcal
F_{q,v}f\|^2_{q,2,v},
$$
and
$$
\|f_t\|^2_{q,2,v}=\frac{1}{t^{2v+2}}\|f\|^2_{q,2,v},\quad\|xf_t\|^2_{q,2,v}=\frac{1}{t^{2v+4}}\|xf\|^2_{q,v,2},
$$
which gives
$$
t^4\|x\mathcal F_{q,v}f\|^2_{q,2,v}+t^2\frac{1}{c^2}\log \left(
B^2_{q,v}\frac{\sigma_2^2}{c^{4(v+1)}}\right)\|f\|^2_{q,v,2}+\|xf\|^2_{q,2,v}\geq0,
$$
and then
$$
\|xf\|_{q,2,v}\|x\mathcal
F_{q,v}f\|_{q,2,v}\geq\psi(c)\|f\|^2_{q,2,v}.
$$
where
$$
\psi(c)=\frac{v+1}{[\sigma_2B_{q,v}]^{\frac{1}{v+1}}}|z_c\log(z_c)|,\quad
z_c=\frac{[\sigma_2B_{q,v}]^{\frac{1}{v+1}}}{c^2},\quad 0<z_c<1.
$$
One can see that
$$
\sup_{0<z_c<1}\psi(c)=\psi(q^\alpha),\quad\alpha=\frac{\log[\sigma_2B_{q,v}]}{2(1+v)\log
q}+\frac{1}{2\log q}.
$$
Let
$$
n_1=\lfloor\alpha\rfloor,\quad n_2=\lceil\alpha\rceil,
$$
where$\lfloor.\rfloor$ and $\lceil.\rceil$ are respectively the
floor and ceiling functions. Now the constant $k$ is given as
follows
$$
k=\psi(q^{n_1}),\quad\text{if}\quad
\lceil\alpha\rceil\geq\alpha-\frac{1}{2\log q}
$$
and
$$
k=\max\{\psi(q^{n_1}),\psi(q^{n_2})\},\quad\text{if}\quad
\lceil\alpha\rceil<\alpha-\frac{1}{2\log q}.
$$
This finish the proof.
\end{proof}

\end{document}